\documentclass[conference]{IEEEtran}
\IEEEoverridecommandlockouts
\usepackage{physics}
\usepackage{cite}
\usepackage{amsmath,amssymb,amsfonts}
\usepackage{algorithmic}
\usepackage{graphicx}
\graphicspath{{pics/}}
\usepackage{textcomp}
\def\BibTeX{{\rm B\kern-.05em{\sc i\kern-.025em b}\kern-.08em
    T\kern-.1667em\lower.7ex\hbox{E}\kern-.125emX}}

\newcommand{\R}{\mathbb{R}}

\newcommand{\Wp}{\mathbb{W}}

\newcommand{\RA}{\rightarrow}
\newcommand{\argmax}{\mathop{\rm argmax}}

\begin{document}

\title{An indirect numerical method for a time-optimal state-constrained control
problem in a steady two-dimensional fluid flow
\thanks{The support from the Russian Foundation for
Basic Research during the projects 16-31-60005, 18-29-03061, and the support of
FCT R\&D Unit SYSTEC -- POCI-01-0145-FEDER-006933/SYSTEC funded by ERDF $|$
COMPETE2020 $|$ FCT/MEC $|$ PT2020 extension to 2018, Project STRIDE
NORTE-01-0145-FEDER-000033 funded by ERDF $|$ NORTE2020, and FCT Project MAGIC
POCI-01-0145-FEDER-032485 funded by ERDF $|$ COMPETE2020 $|$ POCI are
acknowledged. } }

\author{ \IEEEauthorblockN{Roman Chertovskih}
	\IEEEauthorblockA{\textit{Research Center for Systems and Technologies (SYSTEC)}\\
                \textit{Electrical and Computer Engineering Department} \\
		\textit{Faculdade de Engenharia da Universidade do Porto}\\
		Porto, Portugal \\
		roman@fe.up.pt}
	\\
	\IEEEauthorblockN{Nathalie T. Khalil}
	\IEEEauthorblockA{\textit{Research Center for Systems and Technologies (SYSTEC)}\\
                \textit{Electrical and Computer Engineering Department} \\
		\textit{Faculdade de Engenharia da Universidade do Porto}\\
		Porto, Portugal \\
		khalil.t.nathalie@gmail.com}
	\and
		\IEEEauthorblockN{Dmitry Karamzin}
	\IEEEauthorblockA{\textit{Federal Research Center ``Computer Science and Control''} \\
     \textit{of the Russian Academy of Sciences}\\
             Moscow, Russia\\
             dmitry\_karamzin@mail.ru}
	\\
	\IEEEauthorblockN{Fernando Lobo Pereira}
	\IEEEauthorblockA{\textit{Research Center for Systems and Technologies (SYSTEC)}\\
                \textit{Electrical and Computer Engineering Department} \\
		\textit{Faculdade de Engenharia da Universidade do Porto} \\
		Porto, Portugal \\
		flp@fe.up.pt}
}

\maketitle

\begin{abstract}
This article concerns the problem of computing solutions to state-constrained optimal control problems whose trajectory is affected by a flow field. This general mathematical framework is particularly pertinent to the requirements underlying the control of Autonomous Underwater Vehicles in realistic scenarii. The key contribution consists in devising a computational indirect method which becomes effective in the numerical computation of extremals to optimal control problems with state constraints by using the maximum principle in Gamkrelidze's form in which the measure Lagrange multiplier is ensured to be continuous. The specific problem of time-optimal control of an Autonomous Underwater Vehicle in a bounded space set, subject to the effect of a flow field and with bounded actuation, is used to illustrate the proposed approach. The corresponding numerical results are presented and discussed.
\end{abstract}


\section{Introduction}\label{Section_1}

This article addresses the challenges of using computational indirect methods based on the application of the Pontryagin Maximum Principle (PMP) to solve state-constrained optimal control problems arising in the optimal motion planning of Autonomous Underwater Vehicles (AUVs) subject to state and control constraints as well as to the effect of some given fluid vector field. State constraints are specified by the boundary of the free space in which the AUV is allowed to navigate, while the vector flow field, for example, given by currents, moving fronts, and so on, is generated by the diverse underwater phenomena.

In order to facilitate the exposition, we consider a simplified two-dimensional AUV model. That is, the immersion depth is constant, and, moreover, the steering angle can be changed instantaneously. The latter means that the AUV velocity may exhibit discontinuities which can occur at any time due to the control action. This simplified model does not entail any loss of generality of the proposed approach, but it is also still relevant from the point of view of applications. Note, that this approximation is reasonably close  to reality, if the global transition time from a given initial point $A$ to a final point $B$, is several orders of magnitude larger than the time spent on the abrupt change of the vehicle rudder position.

The approach proposed here consists in using numerical algorithms based on an indirect method which employs the PMP to solve state-constrained optimal control problems. As is known, the conventional scheme using this approach encounters a certain computational difficulty in the application of the standard shooting method to solve the two-point boundary value problem which arises from the application of the PMP, due to the Borel measure Lagrange multiplier associated with the state constraints.

How is this difficulty overcome? The key observation consists in the fact that, for the considered AUV model problem and constraints set-up, the regularity condition with respect to the state constraints proposed by R.V. Gamkrelidze in \cite{Gamkrelidze_1959} (see also the classic monograph \cite{Pontryagin_1962}, Chapter 6) is a priori satisfied for any feasible arc. This allows us to assert that the measure Lagrange multiplier from the PMP is continuous, cf. \cite{Arutyunov_Karamzin_2015}. In turn, the regularity condition yields an explicit formula for this multiplier, while, by virtue of its continuity, the junction points at which the optimal arc meets the boundary of the state constraint set, can be found. Thus, largely due to this property, and by using a variation of the standard shooting method (see \cite{nr} for an overview) for numerical solution of the two-point boundary value problem arising from the application of the PMP, an algorithm solving the time-optimal problem is constructed.

There is a vast array of publications on the theory of state-constrained optimal control problems. So far, many important theoretical questions have been investigated. The continuity of the measure-multiplier has also been examined, for example, in \cite{Hager_1979,Maurer_1979,Milyutin_1990,Galbraith_Vinter_2003,Karamzin_Pereira_2019}.
Other contributions made on the general development of this theory can be found,
for example, in \cite{Warga_1964,Dubovitskii_Milyutin_1965,Neustadt_1967,Halkin_1970,Ioffe_Tikhomirov_1979,
Arutyunov_Tynyanskiy_1985,Arutyunov_1991,Vinter_Fereira_1994,Arutyunov_Aseev_1997,
Arutyunov_2000,Vinter_2000,Milyutin_2001,Arutyunov_Karamzin_Pereira_2011,Colombo_2015,
khalil_thesis_2016,Bettiol_Khalil_Vinter_2016}. Issues on numerical solutions to state-constrained problems have been studied, for example, in \cite{Bryson_1969,Betts_1993,Maurer_2000,Haberkorn_2011,Dang_2016,Zeiaee_2017}.
These selective lists of contribution are, obviously, far from exhaustive.

Our article is organized as follows. In Section \ref{Section_2}, the problem formulation for the two-dimensional AUV-motion optimal planning is presented. Section \ref{Section_3} is aimed to discuss the PMP and its application to a particular case of AUV-model problem in which the feasible control set is given by the unit disk in $\R^2$. In Section \ref{Section_4} the numerical results are presented. Section \ref{Section_5} concludes the article with a short summary.

\section{Problem Statement and the Concept of Regularity}\label{Section_2}

Consider an AUV moving in an underwater milieu subject to state constraints and to the influence of fluid flow vector field. The velocity resulting from the force exerted by the fluid flow, represented by the vector $v(x)$ at each point $x$, affect the motion of the vehicle, that is, the way the AUV is propelled.

The problem formulation is described as follows:

\begin{itemize}

\item The waterway is defined on the plane by given affine state constraints.
\item The AUV motion is determined by a linear control system that encompasses the vector field $v(x)$ and the control actuation $u$.
\item The control actuation takes on values in a closed bounded set $U$.
\item The initial and final positions of the AUV are given, respectively, by points $A$ and $B$.
\item The task is to find the minimum time trajectory joining the points $A$ and $B$.
\end{itemize}

More precisely:
\begin{equation}
\begin{aligned}\label{problem: time optimal}
& {\text{Minimize}}
& & T \\
& \text{subject to}
& &  \dot x = u + v(x), \\
&&& x(0)=A,\;\;x(T)=B, \\
&&& -1\le x_1\le 1, \\
&&& u \in U.
\end{aligned}
\end{equation}
Here, $x=(x_1,x_2)$ is the state variable, $u=(u_1,u_2)$ is the control variable. A measurable function $u(\cdot):[0,T]\RA U$, where $U$ is a given compact (the so-called feasible control set), is termed control. The point $A=(a_1,a_2)$ is the starting point, while $B=(b_1,b_2)$ is the terminal point, and $v:\mathbb{R}^2 \to \mathbb{R}^2$ is a smooth map which defines a steady fluid flow varying in space. The terminal time $T$ is free and is supposed to be minimized.

Next, let us formulate the regularity concept w.r.t. the state constraints
for Problem (\ref{problem: time optimal}) (cf. \cite{Karamzin_Pereira_2019}).
Let us define
$$
\Gamma(x,u) = u_1 + v_1(x).
$$
Clearly, $\displaystyle{\frac{\partial \Gamma}{\partial u}(x,u)=(1,0)}$. Then, the regularity condition from \cite{Karamzin_Pereira_2019} is as follows.

\medskip
\noindent
{\bf Assumption R)} Assume that for all $x\in\R^2$ and $u\in U$,
such that $|x_1|=1$, and $\Gamma(x,u)=0$, there exist a vector
$c=(c_1,c_2)\in T_U(u)\cap N^*_U(u)$, and a vector
$d=(d_1,d_2)\in T_U(u)\cap N^*_U(u)$ such that $c_1>0$ and $d_1<0$.

\medskip
Here, $T_U(u)$ is the contingent tangent cone and $N^*_U(u)$ is the
dual of the limiting normal cone. For these definitions, see
\cite{Mordukhovich_2006}. It is clear that, in the case of convex $U$, Assumption R)
always fails to hold wherever there exists a corner point $u_0$ of $U$ for
which $\Gamma(\bar x,u_0)=0$ for some $\bar x=(\bar x_1,\bar x_2)$ such that $|\bar
x_1|=1$.

Nonetheless, Assumption R) represents a kind of a priori verification condition
imposed on the data of the problem. It is easy to single out various
classes of problems which satisfy this condition. One of such class of
problems is demonstrated in the next section. Note that, if Assumption R) is
satisfied, then any feasible arc is regular in the sense proposed in
\cite{Gamkrelidze_1959}, and the measure Lagrange multiplier from the PMP is
continuous. Next, we proceed to the PMP formulation.

\section{Maximum Principle}\label{Section_3}

Consider the extended Hamilton-Pontryagin function
$$
\bar H(x,u,\psi,\mu,\lambda) = \bigl<\psi,u + v(x)\bigr> - \mu \Gamma(x,u) -\lambda,
$$
where $\psi \in \R^2$, $\mu \in \R$ and $\lambda \in \R$.

We assume that Assumption R) holds. Then, for an optimal process $(x^*,u^*,T^*)$, the PMP
ensures the existence of Lagrange multipliers composed by a number $\lambda \in [0, 1]$, an
absolutely continuous adjoint arc $\psi =(\psi_1,\psi_2) \in \Wp_{1,\infty}([0,
T^*]; \R^2)$, and a scalar function $\mu(\cdot)$, such that the following
conditions are satisfied:

\begin{itemize}
\item[(a)]\label{item: adjoint system}Adjoint equation
\begin{align*}
\dot \psi(t) & = -\pdv{\bar H}{x} (x^*(t),u^*(t),\psi(t),\mu(t),\lambda) \\
& =-\psi(t) \frac{\partial v}{\partial x}(x^*(t)) + \mu(t) \frac{\partial v_1}{\partial x}(x^*(t))
\end{align*}
for a.a. $t\in [0,T^*]$;
	
\item[(b)]\label{item: max condition} Maximum condition
\begin{align*}
u^*(t) & \in \argmax_{u\in U} \{{\bar H} (x^*(t),u,\psi(t),\mu(t),\lambda) \} \\
& = \argmax_{u\in U} \{(\psi_1(t)-\mu(t))u_1+\psi_2(t)u_2 \}
\end{align*}
for a.a. $t\in [0,T^*]$;

\item[(c)]\label{item: conservation law} Conservation law
$$
\max_{u\in U} \{ \bar H(x^*(t),u,\psi(t),\mu(t),\lambda) \} =0 \;\; \forall\,t\in [0,T^*];
$$

\item[(d)]\label{item: measure continuity} $\mu(t)$ is constant on the time intervals where
$$
-1 < x_1^*(t) < 1,
$$
increasing on the time intervals where $x_1^*(t)=-1$, and decreasing
on the time intervals where $x_1^*(t)=1$. Moreover, $\mu(\cdot)$
is continuous on $[0,T^*]$;
	
\item[(e)] \label{item: nontriviality condition} Non-triviality condition
$$
\lambda + |\psi_1(t)-\mu(t)| + |\psi_2(t)| > 0\;\;\forall\,t\in [0,T^*].
$$
\end{itemize}

Above, $T^*$ stands for the optimal time, thus, the optimal pair $(x^*,u^*)$ is
considered over time interval $[0,T^*]$.

Now, let us focus on a particular case of the control set $U$. In this article,
we consider the case given by
$$
U:= \mathbb{D}=\{ u =(u_1,u_2) \in \mathbb{R}^2 \ : \ u_1^2+u_2^2 \le 1 \},
$$
that is, by the unit disk in the
plane. The main target is to explicit the necessary optimality conditions
provided by the maximum principle and to determine the appropriate expressions for
the optimal control $u^*$ and the multipliers for a given vector field $v(x)$.
At the same time, it is also required to check whether Assumption R) holds.

After specifying the multipliers, it becomes possible to define an algorithm in order to compute the field of extremals by virtue of the PMP. The
main challenge is that, in the presence of state constraints, the extra
multiplier $\mu(\cdot)$ appears in the PMP, which varies on the subset of $[0,T^*]$ in which the state constraint becomes active, i.e. at all points in time for which the optimal arc reaches the boundary of the state constraint set. If
the regularity condition expressed by Assumption R) does not hold for Problem (\ref{problem: time optimal}), then $\mu(\cdot)$ may have
jumps at such points of time.\footnote{See Example 1 in
\cite{Zakharov_Karamzin_2015}.} Moreover, in the absence of regularity, it is
not clear how to express $\mu(\cdot)$ via the rest of multipliers, that is
$\psi_1(.)$ and $\psi_2(.)$. Both facts entail obvious numerical complexities in
computing the multipliers. Nonetheless, Assumption R) and the continuity of
$\mu(\cdot)$ ensured by this assumption, enable the numerical efficiency of the proposed algorithm, as it will be discussed in Section \ref{Section_4}.

Based on this, consider several claims relevant to the on-going analysis.

\medskip
{\it Claim 1.} Assumption R) holds for Problem (\ref{problem: time optimal})
in which $U=\mathbb{D}$ provided that $|v_1(x)|<1$ $\forall\, x=(x_1,x_2)$ such that $|x_1|=1$.

Indeed, on the boundary of the state constraint one has
$u_1^*(t)=-v_1(x^*(t))$, and hence, by the assumption, it holds that
$|u_1^*(t)|<1$. Then, since
\[T_\mathbb{D}(u) \cap N^*_\mathbb{D}(u)=\{\xi: \bigl<u,\xi\bigr> = 0\} \quad \forall\, u\in\partial\mathbb{D},\]
Assumption R) obviously holds with $d=-c$ with $d_1\ne 0$. The claim is therefore confirmed.

\medskip
Observe that the vehicle actuators should be sufficiently powerful to enable it to cross the given column waterway milieu. This is obviously so, should
the flow field verify the condition $|v_1(x)| < 1$ for all $x$, which implies
the corresponding condition of Claim 1. Moreover, this condition guarantees the existence of a feasible path from $A$ to $B$, as the
main fluid flow is supposed to be along the axis $x_1=0$, in the direction from
$A$ to $B$. Therefore, the application of Filippov's Theorem \cite{Vinter_2000} yields the
following claim.

\medskip {\it Claim 2.}
Problem (\ref{problem: time optimal}) with
$U=\mathbb{D}$ has a solution under the above non-restrictive assumptions
imposed on the vector field $v(x)$.

\medskip
From the PMP, it follows

\medskip
{\it Claim 3.} In the PMP for Problem (\ref{problem: time optimal}), one has
\begin{equation} \label{equation: strict positivity}
|\psi_1(t)-\mu(t)| + |\psi_2(t)| > 0\;\; \forall\,t\in [0,T^*].
\end{equation}

Indeed, if there exists some $t \in [0,T^*]$ such that condition (\ref{equation:
strict positivity}) is violated, then due to the Conservation law (c), we obtain
that $\lambda=0$, which contradicts the Non-triviality condition (e), thus confirming the claim.

\medskip

Due to Claim 2, a solution $(x^*,u^*,T^*)$ to Problem (\ref{problem: time
optimal}) with $U=\mathbb{D}$, exists, while due to Claim 1 the above PMP can be
applied to it. The next step consists in deducing explicit formulas for $u^*$ and $\mu$
expressed via $\psi$. These are needed to solve the four-dimensional boundary
value problem arising from the application of the PMP. From the Maximum condition (b), and by virtue of
(\ref{equation: strict positivity}), the optimal control $u^* = (u_1^*,u_2^*)$
is expressed uniquely via the multipliers (as the control set is strictly
convex) and takes the following form (the dependence on time variable is omitted just to simplify the notation):
\begin{equation}\label{controlo}
u_1^* = \frac{\psi_1 - \mu}{ \sqrt{(\psi_1- \mu)^2 + \psi_2^2}},
\quad   u_2^* = \frac{\psi_2}{\sqrt{(\psi_1 - \mu)^2 + \psi_2^2}}.
\end{equation}

At the boundary points of the state constraint set along the optimal trajectory, we have
$\Gamma(x^*(t), u^*(t)) = 0$, or equivalently, $u_1^*(t) = - v_1(x^*(t))$.
Then,
$$
-v_1^* = \frac{\psi_1 - \mu}{ \sqrt{(\psi_1 - \mu)^2 + \psi_2^2}},
$$
where we set $v_1^*(t):= v_1(x^*(t))$. This implies (bearing in mind that $|v_1^*| < 1$)
\begin{equation}\label{mu}
\mu = \psi_1 + \frac{|\psi_2| v_1^*}{\sqrt{1-v_1^{*2}}}.
\end{equation}
The derived formula holds at the boundary of the state
constraints, that is, on the time intervals on which $|x_1^*(t)|=1$. Moreover,
$\mu(t)$ is increasing on the time intervals where $x_1^*(t)=-1$, decreasing
where $x_1^*(t)=1$, and constant in the interior of the state constraint set $-1< x^*(t) < 1$.

At the same time, $\mu$ can be chosen such that $\mu(0)=0$, (equivalently,
$\mu(T^*)=0$), and it is continuous on $[0,T^*]$. Thus, the above formulas
(\ref{controlo}) and (\ref{mu}) provide  an explicit expression for $u^*$ and $\mu$ via $\psi$, and, thus, also for the boundary value problem
to compute $x^*,\psi$ by the control dynamics of Problem (\ref{problem: time optimal})
and the Adjoint equation (a). This boundary value problem is described and numerically solved in the next section.

\section{The algorithm and numerical results}\label{Section_4}
Summarizing the results discussed above, the field of extremals is
described by the following two-point boundary value problem:
\begin{equation}
\begin{aligned}\label{eqs}
&\dot x = u + v(x), \\
&\dot \psi =-\psi \frac{\partial v}{\partial x}(x) + \mu \frac{\partial v_1}{\partial x}(x),\\
&x(0) =A,\quad x(T)=B, \\
&-1\le x_1\le 1, \\
\end{aligned}
\end{equation}
together with (\ref{controlo}) and (\ref{mu}). Here, $v=(v_1(x_1,x_2),v_2(x_1,x_2))$ is a given vector field satisfying
the condition in Claim 1; points $A=(a_1,a_2)$ and $B=(b_1,b_2)$ are the given initial and final positions, respectively;
travelling time, $T$, is unknown. Note that multiplier $\mu=\mu(t)$ given by relation (\ref{mu}) is continuous and according
to the previous section can be defined in two equivalent ways either $\mu(0)=0$ or $\mu(T)=0$.

\subsection{The algorithm}
The implemented algorithm to solve Problem (\ref{eqs}) numerically is a variation of the shooting method
(see, e.g., \cite{nr} for details) and consists of two parts.

\subsubsection{Backward time integration}
We perform backward time integration by the standard fourth-order Runge-Kutta method
of the ordinary differential equations (ODEs) in (\ref{eqs}) for the initial conditions
in the form\footnote{By virtue of (\ref{equation: strict positivity}),
it is clear to see that one could take $|\psi_{\rm b}(0)|=1$.}:
$$
x_{\rm b}(0)=B,\quad \psi_{\rm b}(0)=(\sin(\theta),\cos(\theta))
$$
for $\theta$ increasing from 0 to $2\pi$ with a step $\delta=0.01$.
Also, $\mu_{\rm b}(t)=0$ for $t\ge0$ while the trajectory stays in the interior of the state constraint set.
For each trajectory $x_{\rm b}(t)$ the distance to the point $A$ is measured as well as the corresponding travelling time.
Local minima in $\theta$ of the distance is found by bisection (i.e. repeating the procedure for a
halved $\delta$ and between the two neighboring values of $\theta$ where the
current estimate of the minimum is computed) such that distance between the trajectory and $A$
is not greater than $10^{-3}$ -- we name such trajectories to be {\it inner extremals}.

If a trajectory meets the boundary, then $\mu_{\rm b}$ is computed on the
boundary using (\ref{mu}). Bisection in $\theta$ is used to find
trajectories minimizing the multiplier at the junction point such that it verifies $|\mu_{\rm b}| < 10^{-3}$.  Along
such trajectories $\mu_{\rm b}(t)$ is continuous and they are potentially parts of
extremals involving a boundary segment -- {\it boundary extremals}.

Trajectories found in this step of the algorithm are shown in
Figs.\ref{f1}-\ref{f3} by red lines.

\subsubsection{Forward time integration}

Next step is forward in time integration of the ODEs in (\ref{eqs})
for the initial conditions in the form:
$$
x_{\rm f}(0)=A,\quad \psi_{\rm f}(0)=(\sin(\theta),\cos(\theta))
$$
also for $\theta\in[0,2\pi)$ with the same step $\delta$. Only trajectories meeting
the boundary with $|\mu_{\rm f}|<10^{-3}$ are of interest. When such
trajectories are computed by bisection in $\theta$, the segment of the trajectory starting from the junction point is computed by integrating the governing equations along the boundary.
Trajectories found at this step of the algorithm are shown in
Figs.\ref{f1}-\ref{f3} by blue lines.

If a trajectory following the boundary meets a point where a backward
in time trajectory meets the boundary with \hbox{$|\mu_{\rm b}|<10^{-3}$},
then the continuity of the multiplier at this time is checked via the following relation:
\[\frac{{\psi_{\rm b}}}{|{\psi_{\rm b}}|} = \frac{{\psi_{\rm f}}-\mu_{\rm f}\begin{bmatrix}
		1 \\
		0
		\end{bmatrix}}{\bigg|\psi_{\rm f}-\mu_{\rm f}\begin{bmatrix}
		1 \\
		0
		\end{bmatrix}\bigg|}. \]
If it is continuous, then the forward in time trajectory
(involving the corresponding boundary segment) and the backward in time
trajectory together constitute an extremal.

\subsection{Numerical results}

\begin{figure}[!]
\centerline{\includegraphics[scale=1.22]{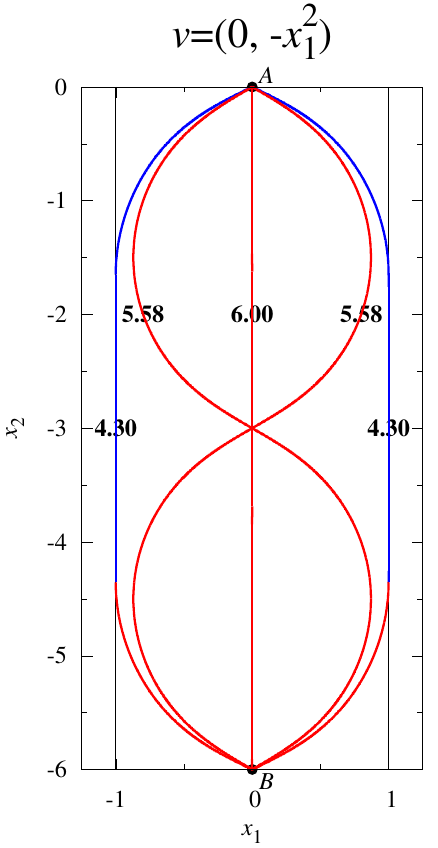}}
\caption{
Field of extremals for $V=(0,-x_1^2)$, $A=(0,0)$ and $B=(0,-6)$.
}
\label{f1}
\end{figure}

\begin{figure}[!]
\centerline{\includegraphics[scale=1.22]{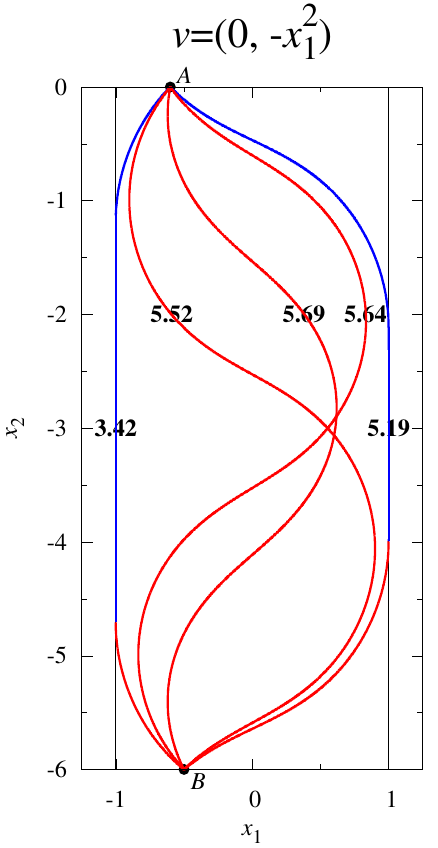}}
\caption{
Field of extremals for $V=(0,-x_1^2)$, $A=(-0.6,0)$ and \hbox{$B=(-0.5,-6)$}.
}
\label{f2}
\end{figure}

\begin{figure}[!]
\centerline{\includegraphics[scale=1.22]{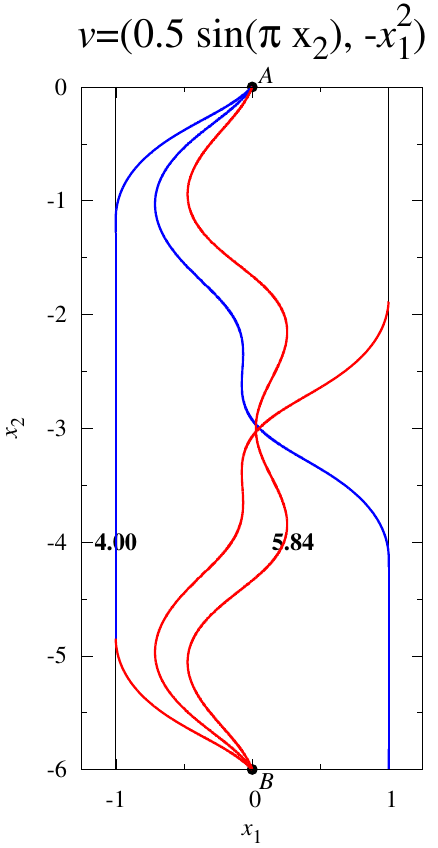}}
\caption{
Field of extremals for $W=(0.5\sin(\pi x_2),-x_1^2)$, $A=(0,0)$ and $B=(0,-6)$.
}
\label{f3}
\end{figure}

For numerical experiments we consider two sample flow velocity fields
\hbox{$V=(0,-x_1^2)$} (Figs.~\ref{f1}-\ref{f2}) and \hbox{$W=(0.5\sin(\pi x_2),-x_1^2)$}
(Fig.\ref{f3}) mimicking real river flows. The former is symmetric
 with vanishing component transversal
to the boundary (permitting validation of the numerical method),
while the latter is more realistic; in both cases, fluid flows are
faster near the boundary.

For the flow $V$, $A=(0,0)$, $B=(0,-6)$, the field of extremals and
the corresponding travelling times are shown in Fig.~\ref{f1} displaying
three inner and two boundary extremals. The flow and the position of the points
$A$ and $B$ are symmetric about the vertical axis
$x_1=0$ and, as a consequence, all the extremals are also symmetric about this
axis. Since the flow is faster at the boundary, as one can expect,
solution to the minimum time problem are the boundary extremals (with optimal traveling
time $T^*=4.3$).

Field of extremals for the same flow, $V$, but for a different position of
the starting and terminal points, $A=(-0.6,0)$ and $B=(-0.5,0)$, is shown
in Fig.~\ref{f2}. Qualitatively, the extremals are of the same nature --
three inner and two boundary extremals, but they are not symmetric anymore.
The left boundary extremal, which is closer
(in comparison to the other boundary extremal) to both points $A$ and $B$, is the solution
to the minimum time problem with $T^*=3.42$. Interestingly, although
the right boundary extremal is more distant from $A$ and $B$ than the inner extremals, it is less time consuming.

For the flow $W$, the field of extremals (see Fig.~\ref{f3}) is constituted of
one inner and one boundary extremals. One backward in time trajectory meeting the
right boundary and one forward in time meeting the same boundary (both are shown
in Fig.~\ref{f3}) do not meet each other and, hence, do not constitute an extremal.
As for the field V above, the minimum time trajectory is the boundary extremal with
$T^*=4$.

\section{Conclusions}\label{Section_5}

A simplified two-dimensional AUV-motion model in an underwater milieu has been considered, albeit under the state constraints. An indirect numerical method based on the application of the Pontryagin maximum principle has been proposed. Formulas for the measure Lagrange multiplier and the extremal control expressed via the state and co-state functions have been derived and the corresponding boundary value problem has been explicated. The numerical results have been presented and discussed.

Besides the presented numerical results concerning the unit disk as a feasible control set, other types of feasible sets have been considered, such as square, ellipsoid, etc. For these models, other types of conditions on the vector field $v(x)$ have been assumed in order to fulfill the regularity requirements.



\end{document}